
\input amstex 
\documentstyle{amsppt}
\input bull-ppt
\keyedby{bull327e/kmt}

\topmatter
\cvol{27}
\cvolyear{1992}
\cmonth{October}
\cyear{1992}
\cvolno{2}
\cpgs{298-303}
\title A shooting approach to the Lorenz 
equations \endtitle
\author S. P.  Hastings and W. C.  Troy \endauthor
\address Department of Mathematics, University of 
Pittsburgh, Pittsburgh,
Pennsylvania 15260\endaddress
\ml SPH\@ MTHSN4.Math.P.H.edu\newline
\indent {\it E-mail address\/}\,: 
Troy\@ VMS.CIS.PITT.EDU\endml
\date January 7, 1992 and, in revised form, April 20, 
1992\enddate
\abstract We announce and 
outline a proof of the existence of a homoclinic orbit of 
the Lorenz 
equations.  In addition, we develop a shooting technique 
and two key 
conditions, which lead to the existence of a one-to-one 
correspondence 
between a set of solutions and the set of all infinite 
sequences of 1's 
and 3's.\endabstract
\endtopmatter

\document


\heading 1.  Introduction\endheading
\par The system of equations discovered by Lorenz {[6]} is 
found 
in computer simulations to have chaotic behavior, by 
practically any 
definition of that term.  A survey appears in [{7}]. 
However, aside 
from local results and various kinds of bifurcation 
analysis, little has 
been proved about these equations. 
\par We have now been able to prove the existence of a 
homoclinic orbit 
for an open set of parameter values. That is, there is a 
nonconstant solution tending to the same equilibrium 
point, in this 
case (0,0,0), at $\pm\infty .$ Such a solution has long 
been conjectured to 
exist and is recognized as an important feature of these 
equations.  
Further, we have proved a theorem, which reduces the 
question of 
whether there are ``chaotic'' orbits to one which, in 
principle, can be 
handled for an open set of parameter values with the 
techniques of 
rigorous numerical analysis, such as interval arithmetic 
[1].  { In other 
words, the amount of computation required for rigorous 
verification 
of the hypotheses of our second theorem is finite; whether 
it is 
practical remains an open question at this time.}  
Computer assisted 
proofs generally leave a gap in understanding, but with 
Theorems 1 
and 2, we believe the gap is smaller than before.  
\par The use of a ``shooting'' technique to obtain a 
homoclinic orbit is 
clearly suggested by numerical integrations. However, its 
analytic 
implementation is difficult and has not been done before, 
to our 
knowledge.  The proof of the second theorem is easier. Its 
conclusion is interesting because it discusses the 
existence of chaos 
without reference to any bifurcation phenomena.  The 
hypotheses 
appear to be true (based on standard numerical 
integrations) for the 
``classical'' parameter values used for the Lorenz 
equations.   
\par Our approach is more elementary than some other 
approaches to 
these equations, such as geometric models, which have 
given deep 
insights about chaotic behavior but have not been shown to 
apply to 
the system which motivated all this work.  Instead of 
Poincar\'e maps 
we use simple one parameter ``shooting.''  Therefore, this 
is a further 
application of methods begun in [5] and continued in [3, 
4].  
Techniques from [8] are also important.
 
\proclaim{Theorem 1} For each $(s,q)$ in some neighborhood 
of the 
point $(10,1)$ there is an $R$ in the interval $(1,1000)$ 
such that the system 
of equations 
$$\aligned x^{\prime}&=s(y-x),\\
y^{\prime}&=Rx-y-xz,\\
z^{\prime}&=xy-qz\endaligned\tag 1$$
has a homoclinic orbit.
 
\endproclaim
 
 We chose $q=1$ as the point to do our analysis instead of 
the 
usual value of 8/3, purely for convenience in the many 
numerical 
calculations.  The lack of precision in the $R$ value 
could easily be 
reduced dramatically with computer assistance.  Also, the 
range of 
parameters $(s,q)$ for which an $R$ exists could be 
expanded greatly;  
but it would be more desirable to find a method that would 
reveal 
the set of such $(s,q)$ analytically, by proving an 
extension of Lemma 2. 
Further homoclinic orbits for a given $(s,q)$ could, in 
principle, 
also be found by our method.  
\heading 2. Outline of proof that there is a homoclinic 
orbit\endheading
\par For any $R>1$, the equilibrium point (0,0,0) is 
unstable and has a 
one-dimensional unstable manifold, which we denote by 
$\gamma$.  We analyze 
the behavior of this manifold for $R$ close to 1 and for 
$R=1000.$ We 
prove two lemmas, which deal with a solution $p(t)=$ 
$(x(t),y(t),z(
t))$ on 
the ``positive'' branch $\gamma^{+}$ of the unstable 
manifold, so that $
p(t)\to (0,0,0)$ 
as $t\to -\infty ,$ and $x(t), y(t)$, and $z(t)$ are all 
positive for large negative  
$t$. 
 \proclaim{Lemma 1} For $(s,q)$ in some open neighborhood 
of $(10,1)$ and 
$R-1$  positive but sufficiently small, $ x(  t),$ $  y
(  t),$ and $  z(  t)$ are all 
positive on the entire line $-\infty <t<\infty .$ 
\endproclaim

\proclaim{Lemma 2} For $(s,q)$ in some open neighborhood 
of $(10,1)$ and 
$R=1000,$ $x(t)$ has at least one zero, $t_1,$ and 
$x^{\prime}$ changes sign 
exactly once in $(-\infty ,t_1].$  
\endproclaim
     

 Lemma 1 is easy to prove analytically, but Lemma 2 is 
much more 
difficult.  With these two results, the existence of a 
homoclinic orbit 
is seen fairly quickly.  We let $R^{*}$ be the infimum of 
all values of 
$R>1$ for which the behavior in Lemma 2 occurs.  For 
$R=R^{*}$ and 
$(x(t),y(t),z(t))$ on $\gamma^{+}$, we can show that 
either (a) $x$ and $
x^{\prime}$ vanish 
simultaneously, or (b) after $\tau_1$, the first zero of 
$x'$, $x$ decreases as long as the 
solution exists but never becomes negative.  This requires 
eliminating 
other options, such as the bifurcation from some finite 
point of new 
zeros of $x^{\prime}$ at $  R^{*}.$  

In case (a), $x=y=0$ at some finite time t, and then the 
uniqueness 
theory for initial value problems implies that $x$ and $y$ 
are identically 
zero.  This is impossible on $\gamma^{+}$.  In case (b), 
it is clear that 
$(x, y, z)$ must approach an equilibrium point as 
$t\to\infty .$ Only a little 
more effort is needed to show that this is the origin and 
the orbit is 
homoclinic.  
 
The proof of Lemma 2 requires a careful study of $\gamma^{+}
.$ The value 
$R=1000$ is, of course, rather arbitrary.  It is not hard 
to prove that 
on some initial interval $(-\infty ,t_0]$, $x,$ $y$, and 
$z$ increase 
monotonically, reaching a point where $y=1$, $0.096\le 
x\le0.1$, and 
$x^2/20<z<0.1$.  From this point the result would follow 
easily with 
computer assistance; but in our opinion, ingenuity, and
considerable effort, is required to follow the solution 
analytically.  
We show that $x$, $y$, and $z$ continue to rise at least 
until they reach 
the levels $z=1000,$ $126.4<x<135.6,$ and $798<y<1000.$ 
From this 
region $y$, and eventually $x$ and $z$, begin to fall.  It 
is apparent that 
at a point where $y$ becomes negative, $z$ must be greater 
than R.  In 
fact, we show that $y=0$ at a point where $155<x<189$ and 
$z>10.4x
.$ 


To obtain these inequalities we use (1) for initial 
estimates 
over suitable intervals, and then iterate to obtain better 
bounds.  Use 
is made of the functions $S=\frac 12(y^2+z^2)-50x^2$ and 
$Q=z-x^2/
20.$ A 
continuation of this process into the region where $y$ is 
negative 
requires some ``tricks'' but finally yields the result.  
The details have 
been submitted elsewhere.


\heading 3.  Criteria for the existence of complicated 
solutions \endheading


 In this section we give a theorem with the conclusion 
that for some 
values of $(s, q, R)$, equation (1) has solutions with 
very complicated 
behavior, in a sense to be made precise.  The hypotheses 
of this 
theorem seem difficult to check analytically; however, the 
result 
seems to us to be an improvement over previous work 
because these 
hypotheses can be confirmed by examining only a compact 
segment of 
$\gamma^{+}$, together with a set of solutions $p$ such 
that $p(0)$ is in a compact 
line segment.  This line segment lies in the intersection 
of the two 
planes $x=y$ and $z=R-1.$ The solutions only have to be 
followed 
over compact time intervals, suggesting that the 
hypotheses can be 
checked rigorously with computer assistance.  


 Standard numerical analysis indicates that the hypotheses 
are 
satisfied, for some parameters, in a robust fashion so 
that the errors 
in the numerical analysis should not be so great that the 
result is 
false.  This reinforces our hope that the theorem can be 
shown to 
apply to (1).  

We have two principle hypotheses for Theorem 2.  The 
first is an extension of Lemma 2.  
 
 \proclaim{Condition A} If p is a solution of $(1)$ with 
$p(0)\in\gamma^{
+},$ then 
$x^{\prime}$ has at least five sign changes and $x$ has at 
least one zero.  If
$\tau_1<\tau_2<\tau_3<\tau_4<\tau_5$ are the first five 
sign changes of $
x^{\prime},$ while $t_1<$ 
$t_2$ are the first two zeros of $x$, then 
$\tau_1<t_1<\tau_2<\tau_
3<\tau_4<\tau_5<t_2.$ 
\RM(If $x$ does not have a second zero, set 
$t_2=\infty$.\RM) \endproclaim\ 
 
This condition is obviously more restrictive than the 
conclusion of 
Lemma 2.  When $(s,q)=(10,1)$, standard numerical 
integrations suggest
that it holds for $R$ approximately in the range 
$(8.2,17.2)$.  If $
q=8/3,$ 
then the $R$ range becomes about (14, 46.6).  
Unpredictable behavior 
exists outside of this range, and a straight forward 
extension of our 
theorem would partly explain this, but be harder to check 
rigorously.  
\par Before stating our second hypothesis, we must 
describe the 
``shooting'' procedure used to obtain the complicated 
solutions.  The 
method is to choose initial conditions $p(0)$ in a certain 
line segment 
in the plane $x=y$ and give an inductive procedure for 
varying $p(0)$
to obtain more and more complex behavior.  
\par To specify this line segment, suppose that Condition 
A is satisfied.  
Then the branch $\gamma^{+}$ of $\gamma$ first crosses the 
plane $
x=y$ at some point 
$p_1,$ which can be shown to lie in the region $z>R-1.$ 
Also, since 
$R>1$, there is an equilibrium point $p_0$ of (1) in this 
plane in the 
positive octant.  At $p_0$, $z=R-1$.  Our shooting set is 
the line 
segment $L$ connecting $p_0$ and $p_1$.  

\par The idea of our second condition is, roughly, that 
solutions starting 
on $L$ do not gain or lose sign changes of $  x^{\prime}$ 
by bifurcation as the 
initial point changes on $L$.  Stated that way, however, 
it appears 
necessary to follow these solutions on $0\le t<\infty ,$ 
clearly not possible 
numerically.  Instead, we consider solutions starting 
along the line 
where zeros of $x^{\prime}$ bifurcate and follow these 
backwards.  This is 
the line M defined by the equations $  x=y$, $  z=R-1$.  
Note that this 
line intersects L only at $  p_0.$ We will explain below 
why this should 
require only a finite amount of computation.  
 
 \proclaim{{\rm Condition B}} Suppose $p$ is a nonconstant 
solution of $(1)$ 
such that $p(0)\in M.$ Then at least one of the 
following is true.  
\roster
\item"(2a)" $p(t)\notin L$ for $t<0$
\item"(2b)" In some interval containing $t=0$, $x\ne 0$ 
and $x^{\prime}$
changes sign four times.
\endroster
 
\endproclaim
     
 To check Condition A numerically, it is first necessary 
to give 
estimates that show that $\gamma^{+}$ intersects a 
specific planar rectangle 
close to, but not including, the origin.  It must then be 
shown that 
every solution starting in this rectangle behaves as 
described in the 
condition.   This process is difficult, because current 
methods of 
interval arithmetic lose about 10 decimal places of 
accuracy for every 
time unit of integration for this system; however it has 
been 
successfully carried out by Hassard and Zhang [2].


To check Condition B, we suggest the use of the well-known 
result 
[7] that the ellipsoid $E$ defined by the inequality 
$$x^2+\frac {10}Ry^2+\frac {10}R(z-2R)^2\le 40R$$
is a positively invariant set for (1) for a range of 
values of $q$ and $s$.  
The line segment $L$ lies in $E$. We consider initial 
points $p(0)$ on the 
(different) line segment $M\cap E$. Assuming that $p$ is 
not constant and 
(2b) cannot be verified, we would integrate from $p(0)$ 
backwards in $t$
and show that the solution leaves $E$ before intersecting 
the line 
segment $L$.  Once $p$ leaves $E$ as $t$ decreases, it 
cannot reenter $E$ at a 
lower $t$ value.  By its nature, a single integration 
using interval 
arithmetic can verify Condition B for an interval of 
initial conditions 
on $M$ around $p(0)$.  The practical difficulty, which up 
to now has 
prevented us from completing this step, is that the length 
of these 
intervals is quite small, so that several thousand initial 
conditions
must be considered. A local analysis around $p_0$ results 
in a 
bound on the length of the time intervals in our 
integrations.  

   We can now state our second theorem.
 
\proclaim{Theorem 2}
Suppose that Conditions {\rm A} and {\rm B} hold for some 
$(s,q,R)$.  
Suppose also that two of the eigenvalues of the linearized 
system around
$p_0$  are complex.  Moreover, 
suppose that $\{M_j\}$ is any infinite sequence of 
$1$\RM's and $3$\RM's. 
Then there 
is a solution $p=(x,y,z)$ of $(1)$ such that $x$ has an 
infinite number of 
zeros in $0<t<\infty ,$ and if $\{t_i\}$ is the sequence 
of consecutive zeros of
$x$ in $[0,\infty )$ and $\sigma_i$ is the number of sign 
changes of $
x^{\prime}$ in $(t_i,t_{i+1}),$ 
then $\sigma_i=M_i$ for $1\le i$ $<\infty$ $.$ \endproclaim
 
\demo{Outline of Proof} Parametrize $L$ by setting 
$p_{\alpha}(0)$ $
=$ 
$\alpha p_0+(1-\alpha )p_1,$ for $0\le\alpha\le 1.$ The 
proof proceeds by induction, choosing 
a sequence of $\alpha$'s giving more and more of the 
prescribed numbers of 
critical points between zeros of $x.$ 
\enddemo

 Because two of the eigenvalues of the linearized system 
around $p_
0$ 
are complex, it follows that if $\alpha$ is close to 1, 
then $p$ crosses the 
plane $y=x$ in $0<t<\infty$ before any possible zero of $ 
x$.  On the other 
hand, for small $\alpha$, $x$ decreases monotonically to 
below 0, after which
$x^{\prime}$ changes sign at least four times before $x=0$ 
a second time.  
Therefore, the first positive zero of $x,$ $t_1(\alpha ),$ 
is defined and 
continuous on some maximal interval of the form 
$[0,\bar{\alpha}$), where $
\bar{\alpha }<1$.  


If $x_{\alpha}$ has at least n positive zeros $t_1(\alpha 
),\dots,t_
n(\alpha ),$ let $t_{n+1}(\alpha )$ denote 
the $(n+1)\roman{st}$ positive zero of $x$ if this exists, 
or else $t_{n+1}
(\alpha )=\infty .$ 
Also, for $1\le i\le n,$ let $\sigma_j=\sigma_j(\alpha )$ 
denote the 
number of sign changes 
of $x^{\prime}$ in $[t_i,t_{i+1}).$ 

\par Suppose that $t_n(\ )$ is continuous on some interval 
$I_n\subset [0
,\bar{\alpha }).$ 
We define three subsets of $I_n$ as follows, where the 
dependence of the $t_j$ and $\sigma_j$ on $\alpha$ is 
again understood:  
$$A_n(I_n)=\{\alpha\in I_n|\text{$t_{n+1}<\infty$,  
$\sigma_n=1$ and $
\sigma_{n+1}\ge 4\}$},$$
$$B_n(I_n)=\{\alpha\in I_n|\text{$t_{n+1}<\infty$,  
$\sigma_n=3$ and $
\sigma_{n+1}\ge 4\}$},$$
$$C_n(I_n)=\{\alpha\in I_n|\text{$t_{n+1}<\infty$,  
$\sigma_n\ge 4$ and $
\sigma_{n+1}\ge 4\}$}.$$
   We prove the following, which imply the theorem. 
\par(i)  Let $I_1=(0,\bar{\alpha })$.  Then $A_1(I_1),$ 
$B_1(I_1)$,  and $
C_1(I_1)$ are all nonempty.
\par(ii)  If, for some n and some $I_n,$ 
$A_n=A_n(I_n),\,\,B_n$, and $C_
n$ are all 
nonempty, then there are intervals $I_{n+1}\subset I_n$ 
and $I^{\prime}_{
n+1}\subset I_n$ such that
\roster
\item"(a)" $t_{n+1}(\ )$ is continuous on $I_{n+1}$ and on 
$I_{n+1}^{
\prime}$;
\item"(b)" $I_{n+1}\cap A_n$ and $I_{n+1}^{\prime}\cap 
B_n$ are nonempty;
\item"(c)" The sets $A_{n+1}(I_{n+1})$, $B_{n+1}(I_{n+
1})$, $C_{n+1}(
I_{n+1})$, $A_{n+1}(I_{n+1}^{\prime}),B_{n+1}(I_{n+
1}^{\prime})$,\newline
and $C_{n+1}(I_{n+1}^{\prime})$ are all non-empty.
\endroster

 We do not have space for the details here, and they have 
been
submitted elsewhere.   Condition B is used to show that as 
$\alpha$ 
varies, the number of sign changes of 
$x_{\alpha}^{\prime}$ between consecutive zeros 
of  $x$ can decrease from four or more to three, or from 
two or three 
to one, only when the fourth or second of these sign 
changes tends 
to infinity on the $t$ axis.  The number of sign changes 
of $x^{\prime}$ cannot 
jump directly from four or more to one without passing 
through an 
open set of $\alpha$'s where there are three.  

These theorems have a few simple corollaries, which we will 
mention when details of the proofs are published.  Since 
the 
divergence of the Lorenz vector field is negative, volumes 
are 
reduced by the flow.  We hope to investigate whether our 
results
have any consequences about the existence of ``strange 
attractors.''
 


\enddocument